\documentclass[12pt]{article} \usepackage{amstext,amsfonts,graphics}
\usepackage{amsmath,amssymb,graphicx,latexsym}
\newtheorem{thm}{Theorem} \newtheorem{lemma}[thm]{Lemma}
 
\newenvironment{pf} {\noindent{\sc Proof. }}{{\hfill
$\Box$}\par\vskip2\parsep} 
\newenvironment{pfof}[1]
{\par\vskip2\parsep\noindent{\sc Proof of\ #1. }}{{\hfill $\Box$}
  \par\vskip2\parsep}

\newcommand{\py}{\mbox{\bf P}\,}

\newcommand{\z}{{\mathbb Z}}

\newcommand{\dof}{\bf\boldmath}

\newcommand{\nat}{\mathbb{N}}
\newcounter{mycount}
\newenvironment{mylist}{\begin{list}{(\roman{mycount})}%
{\usecounter{mycount}\itemsep 0pt}}{\end{list}}

\title{Universal finitary codes with exponential tails}
\date{\today}

\author{Nate Harvey
\and Alexander E. Holroyd\thanks{
Research funded in part by NSF (USA) grants DMS-0072398, CCR-0121555, NSERC 
(Canada) grant, CPAM, MSRI and BIRS.}
\and Yuval Peres\thanks{
Research supported in part by NSF Grants
DMS-0104073 and DMS-0244479.}
\and Dan Romik
}

\begin{document}

\maketitle
\renewcommand{\thefootnote}{}

\begin{abstract}
    In 1977, Keane and Smorodinsky \cite{ks:hom} showed that there exists a
    finitary homomorphism from any finite-alphabet Bernoulli process to any
    other finite-alphabet Bernoulli process of strictly lower entropy.  In
    1996, Serafin \cite{serafin} proved the existence of a finitary
    homomorphism with finite expected coding length.  In this paper, we
    construct such a homomorphism in which the coding length has
    exponential tails. Our construction is source-universal, in the sense
    that it does not use any information on the source distribution other
    than the alphabet size and a bound on the entropy gap between the
    source and target distributions.  We also indicate how our methods can
    be extended to 
    prove a source-specific
    version of the result for Markov chains.
\end{abstract}

\section{Introduction}
\label{intro}

Let $a, b\in\nat$, and define the two finite alphabets ${\bf A}=\{i \in \z
: 1 \leq i \leq a\}$, ${\bf B}=\{i \in \z : 1 \leq i \leq b\}$. Equip the
sequence spaces ${\bf A}^\z$ and ${\bf B}^\z$ with the product
$\sigma$-algebras ${\cal A}$ and ${\cal B}$ respectively. A measurable map
$\varphi:\Omega\to {\bf B}^\z$, where $\Omega\subseteq {\bf A}^\z$ is
measurable, is called {\dof translation-equivariant} if for all
$x=(x_i)_{i\in\z} \in\Omega$, the left shift ${\cal
   T}(x)=(x_{i+1})_{i\in\z}$ of $x$ is also in $\Omega$ and the equality
$\varphi({\cal T}(x))={\cal T}(\varphi(x))$ holds. A
translation-equivariant map $\varphi:\Omega\to{\bf B}^\z$ is {\dof
   finitary} if for all $x\in\Omega$, there exists an $N\in\nat$ such that
for all $y\in\Omega$, if $(x_i)_{|i|\le N} = (y_i)_{|i|\le N}$ then
$\varphi(x)_0 = \varphi(y)_0$.  In this case we let $N_\varphi(x)$ be the
minimal such $N$, and call $N_\varphi$ the {\dof coding length} of
$\varphi$.

If $p=(p(i))_{i\in{\bf A}}$ is a probability vector (that is to say,
$p(i)\ge 0$ for all $i\in{\bf A}$, and $\sum_{i\in{\bf A}} p(i)=1$), let
${\bf P}_p$ be the product measure $p^\z$ on ${\cal A}$. The quadruple
$B(p)=({\bf A}^\z, {\cal A}, {\bf P}_p, {\cal T})$ is the {\bf Bernoulli
   shift} of $p$. Similarly if $q=(q(i))_{i\in{\bf B}}$ is a probability
vector on ${\bf B}$, let ${\bf P}_q$ be the product measure $q^\z$ on
${\cal B}$, and let $B(q)=({\bf B}^\z, {\cal B}, {\bf P}_q, {\cal T})$ be
the Bernoulli shift of $q$. A {\dof homomorphism} $\varphi$ from $B(p)$ to
$B(q)$ is a translation-equivariant map $\varphi:\Omega\to {\bf B}^\z$ 
with the properties that
$\Omega\in {\cal A}$, ${\bf P}_p(\Omega)=1$, and for all $E\in{\cal
   B}$ we have ${\bf P}_p(\varphi^{-1}(E)) = {\bf P}_q(E)$.

Denote by $h(p)=-\sum_i p(i)\log p(i)$ the entropy of a probability vector.
Keane and Smorodinsky \cite{ks:hom} proved that if $h(p)>h(q)$, then there
exists a finitary homomorphism $\varphi$ from $B(p)$ to $B(q)$.  Serafin
\cite{serafin} demonstrated that $\varphi$ may be chosen in such a way that
the expected coding length ${\bf E}_p(N_\varphi)$ is finite. Iwanik and
Serafin \cite{iwanser} strengthened this result to all moments below the
second.

We say that a finitary homomorphism has {\dof exponential tails} if
there exist $c>0$ and $0< d < 1$ such that ${\bf P}_p(N_\varphi \geq n) \leq
c\cdot d^n$ for all $n$.  In general, we say that a non-negative
sequence $(c_n)_{n \in {\mathbb N}}$ {\dof decays exponentially} if
there exist $c>0$ and $0<d<1$ such that $c_n \leq c\cdot d^n$ for all
$n$.  We say that a random variable $W$ has {\dof
exponential tails} if $(\py(|W| \geq n))_{n \in {\mathbb N}}$ decays
exponentially.

Our main result is a new construction of a finitary homomorphism from
$B(p)$ to $B(q)$, when $h(p)>h(q)$. Our construction improves on the
above-mentioned results in two ways. First, the coding length will have
exponential tails. Second, the homomorphism is source-universal, in the
sense that the same function works simultaneously for all source vectors
$p$ over a given alphabet which have full support and whose entropy is
greater than $h(q)$ by at least a given $\varepsilon$. In particular, this
answers the open problem mentioned in the last two lines of \cite{peres}.

The precise result is the following:

\begin{thm}
\label{main}
Fix a probability vector $q=(q(i))_{i\in{\bf B}}$, and fix $\varepsilon>0$.
There exists a measurable subset $\Omega\subseteq A^\z$ and a finitary
translation-equivariant map $\varphi:\Omega\to {\bf B}^\z$, such that for
any probability vector $p=(p(i))_{i\in{\bf A}}$ for which $p(i)>0$ for all
$i\in{\bf A}$ and $h(p)\ge h(q)+\varepsilon$, $\varphi$ is a finitary
homomorphism from $B(p)$ to $B(q)$ with exponential tails.
\end{thm}

\bigskip
Here is a brief description of the motivation and method of the proof of
Theorem \ref{main}. A homomorphism can be thought of as a
translation-equivariant function that, given a sequence of independent
samples with distribution $p$, simulates a sequence of independent
samples with distribution $q$. Thus, it is natural to try to construct
such a function using existing constructions of functions that simulate
one discrete distribution using another. Such constructions have been
described in \cite{elias}, \cite{hanhoshi}, \cite{yao}, \cite{romik},
\cite{peres}.

Our construction combines elements from several of these constructions.
First, divide the source sequence $(x_k)_{k\in\mathbb{Z}}$ into
{\bf blocks} separated by {\bf markers}.  A marker is defined as an
appearance of
a certain (sufficiently rare) pattern, say a 2 followed by $t$ 1's, where
$t$ is a large enough integer. Next, feed the contents of each block
into a specially designed function which converts these 
approximately independent $p$-distributed samples into independent
unbiased bits. Next, at each block, start feeding these unbiased bits
into another function designed to generate a number of independent
samples of the distribution $q$ sufficient to fill the length of that
block. This function may require more bits than that block contains,
but on the average, because of entropy considerations it requires less.
Any unused bits are then used to satisfy blocks whose simulation did
not end in the first round. Continuing in this manner one
obtains the required number of samples of $q$, which are then used to
generate the value $\varphi(x)$. Everything is done simultaneously for all blocks
in a
translation equivariant manner, with an added bonus being the source
universality property.

\bigskip

The following is an extension of Theorem \ref{main} to Markov chains.

\begin{thm}
\label{markov}
Let $\alpha=(\alpha_{i,j})_{i,j\in{\bf A}}$,
$\beta=(\beta_{i,j})_{i,j\in{\bf B}}$ be two aperiodic, irreducible Markov
transition matrices over the finite alphabets ${\bf A}, {\bf B}$. Let
${\cal M}(\alpha)=({\bf A}^\z, {\cal A}, P_{\alpha}, {\cal T})$, ${\cal
   M}(\beta)=({\bf B}^\z, {\cal B}, P_{\beta}, {\cal T})$ be the stationary
Markov shifts of $\alpha$ and $\beta$ respectively, and denote their
entropies by $h(\alpha), h(\beta)$. If $h(\alpha)>h(\beta)$, then there
exists a finitary homomorphism from ${\cal M}(\alpha)$ to ${\cal M}(\beta)$
with exponential tails.
\end{thm}

We indicate in Section 5 how our methods may be adapted to prove Theorem
\ref{markov}.
For Markov chains, our construction is not source-universal, except
in the weak sense that it will work, under the assumption of an
entropy gap, simultaneously for all Markov transition matrices with all
entries positive -- see Section 5.

%
%
%
%
%

\section{Simulations} \label{simulations}

In this section, we construct two procedures for simulating one
discrete distribution from another. The constructions are variants of
those used by Elias \cite{elias}, Han and Hoshi \cite{hanhoshi}, Knuth
and Yao \cite{yao}, and Romik \cite{romik}. In all of these
constructions, a key point is that the loss in entropy is small.

\subsection{Simulating a distribution from independent unbiased random
bits}

For $b \in \nat$, let ${\bf B}= \{ i \in {\z}:1 \leq i \leq b \}$ be
a finite alphabet, and let $q=(q(i))_{i \in {\bf B}}$ be a probability
vector. Let $(\{0,1\}^{\nat}, {\cal F}, {\dof P})$ be the probability
space of (one-sided) infinite binary strings, equipped with the
natural product $\sigma$-algebra, and the probability measure under
which the coordinate functions are independent unbiased random
bits. Let ${\dof E}$ denote expectation with respect to the measure
${\dof P}$.

A {\dof simulation of $q$ from independent unbiased bits} is a pair of
measurable functions $T:\{0,1\}^{\nat}\to \nat$ and $S:\{0,1\}^{\nat}\to
{\bf B}$, defined ${\dof P}$-a.s., with the following properties:
\begin{mylist}
\item If $x=(x_i)_{i\in\nat}$ and $\tilde{x}=(\tilde{x}_i)_{i\in\nat}$
are elements of $\{0,1\}^{\nat}$ such that $(\tilde{x}_1, \ldots,
\tilde{x}_{T(x)}) = (x_1,\ldots,x_{T(x)})$, then $T(\tilde{x})=T(x)$
and $S(\tilde{x})=S(x)$.
\item Under the measure ${\bf P}$, $S(x)$ has distribution $q$.
\end{mylist}
$T$ is called the {\dof stopping time} of the simulation, and $S$ is
called the {\dof output symbol}.

The following theorem was first proved by Knuth and Yao
\cite{yao}. We present an independent proof with a more explicit
construction.

\begin{thm} \label{simu1}
There exists a simulation $(T,S)$ of $q$ from independent unbiased bits
satisfying the additional properties:
\begin{mylist}
\item $T(x)$ has exponential tails. More precisely,
${\bf P}(T(x)>k) \le \frac{b+1}{2^k}$.
\item ${\bf E}(T(x)) \le \frac{h(q)}{\log 2} + 6$.
\end{mylist}
\end{thm}

\begin{pf}
Construct $T$ and $S$ as follows. Define a partition of $[0,1]$ by $0 =
Q_0 < Q_1 < \ldots < Q_b = 1$, where $Q_j = \sum_{i=1}^j
q(i)$. Define
$$ T(x) = \min \left\{ k\in \nat : \textrm{for some }1\le j \le b,\ \,
Q_{j-1} < \sum_{i=1}^k \frac{x_i}{2^i} < 
Q_j - \frac{1}{2^k} \right\}, $$
$$ S(x) = \textrm{the unique }1\le j\le b\textrm{ for which }
Q_{j-1} < \sum_{i=1}^{T(x)} \frac{x_i}{2^i} < Q_j -
\frac{1}{2^{T(x)}}. $$ 
In words, the idea is to consider $x=(x_i)_{i\in\nat}$ as the binary
expansion of a number $u=\sum_{i=1}^\infty x_i 2^{-i} \in [0,1]$, and
to define the output symbol $S(x)$ as that $1\le j\le b$ for which
$u\in(Q_{j-1}, Q_j)$. Determining the correct $j$ necessitates
looking at only the first $T(x)$ bits in the binary expansion of $u$.
So $T(x), S(x)$ are defined for all $x$ which are not the binary
expansions of any of the $Q_j$, and property (i) in the
definition of a simulation is clearly satisfied.  Also, since,
under the measure ${\dof P}$, $u$ is uniformly distributed in $[0,1]$,
we have that
$$ {\dof P}(S(x) = j) = \textrm{Lebesgue measure of }(Q_{j-1},
Q_j) = q(j), $$ 
so property (ii) is also satisfied, and $(T,S)$ is indeed a simulation
of $q$ from independent unbiased bits. To prove the additional
properties claimed in the theorem, note that
$$
{\dof P}(T(x)>k) = \py\bigcup_{j=0}^b \left\{Q_j \in
\Big(\sum_{i=1}^k \frac{x_i}{2^i},\sum_{i=1}^k
\frac{x_i}{2^i}+\frac{1}{2^k}\Big) \right\}  \le \frac{b+1}{2^k},
$$
establishing \ref{simu1}(i). For \ref{simu1}(ii), let for $0\le j\le b$
$$ Q_j = \sum_{i=1}^\infty a_{j,i} 2^{-i}, \qquad (a_{j,i} \in
\{0,1\}, i\in \nat )$$
be a binary expansion of $Q_j$, and define for $1\le j\le b$
$$ m_j = \min\{ k\in \nat : a_{j-1,k} \neq a_{j,k} \}
= \left\lceil -\frac{\log (Q_j - Q_{j-1})}{\log 2} \right\rceil. $$
Then, checking the definitions we see that, for $l\ge m_j$,
\begin{multline}\label{eq:event}
\{ x : S(x)=j, T(x)=l \} \\
= \{ x : (x_1,\ldots,x_l) = (a_{j-1,1},\ldots,a_{j-1,l-1},1),
a_{j-1,l}=0 \} \\ \cup\ 
\{ x : (x_1,\ldots,x_l) = (a_{j,1},\ldots,a_{j,l-1},0),
a_{j,l}=1 \},\qquad\,
\end{multline}
while for $l< m_j$, this event is empty. Since
$$ {\bf E}(T(x)) = \sum_{j=1}^b \sum_{l=1}^\infty l\cdot {\bf P}(T(x)=l,
S(x)=j), $$ 
\ref{simu1}(ii) will follow if we prove that for all $1 \le j\le b$,
$$ 
\sum_{l=1}^\infty l\cdot {\bf P}(T(x)=l, S(x)=j) \le
-\frac{q(j)\log q(j)}{\log 2} + 6 q(j).
$$
Indeed, by using \eqref{eq:event} we see that the representation
$$ q(j) = {\bf P}(S(x) = j) = \sum_{l=m_j}^\infty
{\bf P}(S(x)=j, T(x)=l) $$
can be rewritten as a representation of $q(j)$ as a sum of negative
powers of $2$, namely
$$ q(j) = \sum_{i=1}^\infty 2^{-n_{j,i}}, $$
where each summand $2^{-n_{j,i}}$ is the probability of one of the
two events on the right-hand side of \eqref{eq:event}. Arrange
the $n_{j,i}$ such that $n_{j,1} \le n_{j,2} \le n_{j,3} \le ... $, and
note also that $n_{j,i} < n_{j,i+2}$ for all $i \in \mathbb{N}$, since
any given power of $2$ appears at most twice in the sum.
Then in particular we get
\begin{eqnarray*}
2^{-n_{j,1}}\ \, \le\ \, 
q(j) &\le& 2^{-n_{j,1}} + 2^{-n_{j,1}} + 2^{-n_{j,1}-1} + 2^{-n_{j,1}-1}
+ 2^{-n_{j,1}-2} + \ldots \\
&=& 2 \left(1+\frac{1}{2}+\frac{1}{4}+\ldots\right)2^{-n_{j,1}} =
2^{2-n_{j,1}},
\end{eqnarray*}
so $n_{j,1} < -\frac{\log q(j)}{\log 2}+2$. This then gives
\begin{eqnarray*}
\sum_{l=1}^\infty l\cdot {\bf P}(T(x)=l, S(x)=j) &=&
\sum_{i=1}^\infty n_{j,i}2^{-n_{j,i}} \qquad\qquad\qquad\qquad
\end{eqnarray*}

\vspace{-22.0pt}
\begin{eqnarray*} \qquad &=&
\sum_{i=1}^\infty n_{j,1}2^{-n_{j,i}} +
\sum_{i=1}^\infty (n_{j,i}-n_{j,1})2^{-n_{j,i}} \\ &\le&
\left(-\frac{\log q(j)}{\log 2}+2\right) \sum_{i=1}^\infty
2^{-n_{j,i}} + 
2^{-n_{j,1}}\cdot 2\sum_{k=0}^\infty k 2^{-k}
\\ &\le&
-\frac{q(j)\log q(j)}{\log 2} + 2 q(j)+4 q(j),
\end{eqnarray*}
as required.
\end{pf}

\paragraph{Remarks.} Knuth and Yao \cite{yao} proved Theorem
\ref{simu1} with the better constant 2 replacing 6 in Theorem
\ref{simu1}(ii). Their construction is described in slightly less
concrete terms than the one above, but it is optimal, in the sense
that the constant 2 is sharp, and in the strong sense that the
stopping time is stochastically dominated by the stopping time of
\emph{any} possible simulation of $q$ using independent unbiased
bits. For more information see Section 5.12 of \cite{coverthomas}.

A generalization of the construction given above
to simulation of $q$ from a source of independent
samples of an arbitrary probability vector $p$ was given by Han and
Hoshi \cite{hanhoshi} and independently by Romik \cite{romik}, both of
whom proved a statement analogous to Theorem \ref{simu1}(ii), with the
base-2 entropy of $q$ replaced by the ratio $h(q)/h(p)$, and with the
constant 6 replaced by some function of the vector $p$.

\subsection{Simulating independent unbiased random bits from a
block with an excluded pattern}

For $a \in \nat$, let ${\bf A}= \{ i \in {\z}:1 \leq i \leq a \}$ be
a finite alphabet, and let $p=(p(i))_{i \in {\bf A}}$ be a probability
vector. For each $n\in\nat$, let $({\bf A}^n, p^n)$ be the discrete
elementary probability space of ${\bf A}$-valued $n$-tuples, with the
i.i.d. product measure with marginal probabilities $p$. In the case $a=2$ of a binary distribution,
Elias \cite{elias} constructed a function that, given an ${\bf
A}^n$-valued, $p^n$-distributed input, simulates a \emph{random
number} of independent unbiased random bits; that is, a pair $(N,F)$,
where $N$ is an $\nat$-valued random variable, and for each $k\in\nat$,
given $N=k$, the random vector
$F$ is distributed uniformly on the set $\{0,1\}^k$.

We shall need a generalization of Elias's construction, which takes as
input $n$ independent samples from a general discrete distribution,
conditioned on the non-appearance of a certain pattern, and returns a
random number of independent unbiased random bits.

For any $t\in\nat$, let $E_{n,t}$ be the subset of ${\bf A}^n$
consisting of all vectors $x=(x_1,\ldots,x_n) \in {\bf A}^n$ for which
for no $i\in\{1,\ldots,n-t\}$ is it true that
$$ x_i = 2,\  x_{i+1}=x_{i+2}=\ldots=x_{i+t-1}=1. $$
That is, $E_{n,t}$ contains all vectors which, considered as words, do
not contain the pattern ``2 followed by $t-1$ 1's''. We sometimes
call such vectors ``pattern-free''. Let
$\tilde{p}_{n,t}$ be the measure $p^n$ conditioned on $E_{n,t}$.
Let $\{0,1\}^* = \cup_{k=0}^\infty \{0,1\}^k$ be the set of finite
strings over the alphabet $\{0,1\}$.

\begin{thm}\label{simu2}
For any $n, t\in \mathbb{N}$
there exist functions 
$$ \begin{array}{l}N_{n,t}:E_{n,t}\to\{0,1,2,\ldots\}, \\
F_{n,t}:E_{n,t}\to \{0,1\}^*, \\
G_{n,t}:E_{n,t}\to \{1,2,\ldots,n^{a-1}\}\end{array}$$
with the properties:
\begin{mylist}
\item Under the measure $\tilde{p}_{n,t}$, for each $k\in
\{0,1,2,\ldots\}$ for which\\ $\tilde{p}_{n,t}(N_{n,t}(x)=k)>0$,
conditioned on $N_{n,t}(x)=k$, the random vector
$F_{n,t}(x)$ is uniformly distributed
on $\{0,1\}^k$.
\item The function $x\to(N_{n,t}(x),F_{n,t}(x),G_{n,t}(x))$ is injective.
\item For any $x\in E_{n,t}$ we have $N_{n,t}(x) \le (\log a/\log 2) n$.
\end{mylist}
\end{thm}

\bigskip
Before going on with the proof, we explain briefly the idea behind this
construction and its importance in what follows. 
The functions $N_{n,t},F_{n,t},G_{n,t}$ accept as input
a $\tilde{p}_{n,t}$-distributed random variable and produce a binary
string, $F_{n,t}(x)$, of length $N_{n,t}(x)$. Conditioned on the number
of bits, the binary string is distributed uniformly over all binary
strings of that length, in other words contains $N_{n,t}(x)$ independent
unbiased bits. We would like to ensure that the construction is efficient,
i.e., extracts enough information from the input; this is guaranteed by
claim (ii), which
states that adding the complementary information $G_{n,t}(x)$
makes the function injective, and so enables reconstruction of
the input, together with the statement that the range of $G_{n,t}(x)$
is relatively small, so the amount of entropy contained in it is limited.
As for claim (iii), it will be used in our proof that the homomorphism
we construct has exponential tails.

\bigskip

\begin{pfof}{Theorem \ref{simu2}}
Throughout the proof, for convenience we shall consider $n$
and $t$ as fixed and in most places omit reference to the dependence
of the various quantities on them. To construct the functions
$N_{n,t}, F_{n,t}, G_{n,t}$, we first divide $E_{n,t}$ into classes of
equiprobable elements. We do this as follows. Let
$$ C = \{ m=(m_1,m_2,\ldots,m_a)\in \z^a : m_i\ge 0,\ 1\le i\le a,\ \,
m_1+\ldots+m_a = n \}. $$
For $x=(x_1,\ldots,x_n)\in E_{n,t}$ and $1\le i\le a$, let
$$ c_i(x) = \# \{ 1\le j\le n : x_j = i \}, $$
and let
$$ \textrm{count}(x) = (c_1(x),\ldots,c_a(x)) \in C. $$
Then clearly $E_{n,t}$ can be written as the disjoint union
$$ E_{n,t} = \bigcup_{m\in C} \{ x\in E_{n,t} : \textrm{count}(x)=m
\} =: \bigcup_{m\in C} D_m. $$
For each $m\in C$, we have for all $x\in D_m$ that
$$ p^n(x) = p(1)^{m_1} p(2)^{m_2} \ldots p(a)^{m_a}, $$
so
$$ \tilde{p}_{n,t}(x) = \frac{p(1)^{m_1} p(2)^{m_2} \ldots p(a)^{m_a}}
{p^n(E_{n,t})}. $$
In other words, all elements of $D_m$ are equiprobable under
$\tilde{p}_{n,t}$. Now, for each $m\in C$, let
$d_m=|D_m|$, and write
$$ d_m = \sum_{i=1}^{s_m} 2^{r_{m,i}},
\qquad r_{m,1} > r_{m,2} > \ldots > r_{m,s_m}\ge 0,
$$
for the binary expansion of $d_m$. The functions $N_{n,t}$, $F_{n,t}$
and $G_{n,t}$ may now be defined as follows. For each $m\in C$,
arrange the elements of $D_m$ in lexicographical order, and for each
$x\in D_m$ denote by $\textrm{rank}(x)$ the position of $x$ in this order.
Set for each $x\in D_m$
\begin{eqnarray*}
N_{n,t}(x) &=& r_{m, k^*(x)}, \quad
k^*(x) = \min\left\{ 1\le k\le s_m : \sum_{i=1}^k 2^{r_{m,i}} \ge
 \textrm{rank}(x) \right\}, \\
F_{n,t}(x) &=& \textrm{the length-$N_{n,t}(x)$ binary expansion of the
number }\\ & &
\qquad \sum_{i=1}^{k^*(x)} 2^{r_{m,i}}-\textrm{rank}(x), \\
G_{n,t}(x) &=& \textrm{the position
of $m$ in the lexicographical order on $C$}.
\end{eqnarray*}
The functions $N_{n,t}$, $F_{n,t}$, $G_{n,t}$ are clearly defined on
all $E_{n,t}$ and have the desired range. In words, we have used the
lexicographical order to give an explicit partition of $D_m$ into
subsets of sizes $2^{r_{m,i}},\ \ 1\le i\le s_m$. On each subset of
size $2^{r_{m,i}}$ we define $N_{n,t}=r_{m,i}$, and for the value of
$F_{n,t}$ assign to the $2^{r_{m,i}}$ possible elements the
$2^{r_{m,i}}$ different binary strings of length $r_{m,i}$. 
This implies claim \ref{simu2}(i), since the elements of $D_m$ are
equiprobable. The
function $G_{n,t}$ is defined so as to encode the residual information
needed to recover the value of $x$ given $F_{n,t}(x)$. Indeed, if $x, x' \in
E_{n,t}$ and $x\neq x'$, then either $x\in D_m$, $x\in D_{m'}$ for
some $m\neq m'$, in which case clearly $G_{n,t}(x)\neq G_{n,t}(x')$,
or $x,x'$ are in the same $D_m$ but $\textrm{rank}(x)\neq
\textrm{rank}(x')$, whence $F_{n,t}(x)\neq F_{n,t}(x')$. This proves
claim \ref{simu2}(ii). Finally, claim \ref{simu2}(iii) is immediate, since
for all $m\in C$, $1\le i\le s_m$ we have
$2^{r_{m,i}} \le d_m \le |E_{n,t}| \le a^n$, so
$N_{n,t}(x)=r_{m,k^*(x)} \le (\log a/\log 2)n$. 
%
%
%
%
%
%
%
%
%
%
%
\end{pfof}

\section{Construction of the homomorphism}

\begin{figure}
\begin{center}
\begin{picture}(300,370)(0,0)
\put(30,360){\ldots 21113132232111211133312111111122213312111\ldots}
\put(80,350){\vector(0,-1){30}}
\put(100,330){Division into blocks}
\put(0,300){\ldots\ \ 2111\ 313223\!\ \ \ 2111\ \ \ \ 2111\ 3331\ \ \ \ 2111\
11112221331\ \ \ \ 2111\ \ \ldots}
\put(20,297){\framebox(28,13)}
\put(94,297){\framebox(28,13)}
\put(133,297){\framebox(28,13)}
\put(200,297){\framebox(28,13)}
\put(308,297){\framebox(28,13)}
\put(22,288){$\times$\scriptsize=marker}
\put(96,288){$\times$}
\put(135,288){$\times$}
\put(202,288){$\times$}
\put(310,288){$\times$}
\put(80,280){\vector(0,-1){30}}
\put(100,260){Computation of the associated bit strings}
\put(22,220){$\times$}
\put(96,220){$\times$}
\put(135,220){$\times$}
\put(202,220){$\times$}
\put(310,220){$\times$}
\put(0,232){\ldots\ \ 0010111\ \ \ \ \ \ \ \ \,$\phi$ \ \ \ \ \ \ \
$\phi$\ \ \ \ \ \ \ \ \ \ \ \ \ \ \ \!\ 1011010001\ \ \ \ \ \ \ \ \ \ \
\ 11\ \ \ldots}
\put(80,210){\vector(0,-1){30}}
\put(100,190){Step $(3,0)$: The simulators are running}
\put(22,130){$\times$}
\put(96,130){$\times$}
\put(135,130){$\times$}
\put(202,130){$\times$}
\put(310,130){$\times$}
\put(0,142){\ldots\ \ 0010111\ \ \ \ \ \ \ \ \,$\phi$ \ \ \ \ \ \ \
$\phi$\ \ \ \ \ \ \ \ \ \ \ \ \ \ \ \!\ 1011010001\ \ \ \ \ \ \ \ \ \ \
\ 11\ \ \ldots}
\put(20,150){\includegraphics{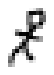}}
\put(18,166){\tiny 1}
\put(32,160){\vector(1,0){20}}
\put(175,155){\includegraphics{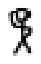}}
\put(172,170){\tiny 2}
\put(186,155){\includegraphics{queueup.eps}}
\put(185,170){\tiny 3}
\put(200,150){\includegraphics{sticks.eps}}
\put(200,169){\tiny 4}
\put(212,160){\vector(1,0){20}}
\put(308,150){\includegraphics{sticks.eps}}
\put(308,169){\tiny 5}
\put(320,160){\vector(1,0){20}}
\put(80,120){\vector(0,-1){30}}
\put(100,100){...\ Step $(3,8)$\ ...}
\put(22,50){$\times$}
\put(96,50){$\times$}
\put(135,50){$\times$}
\put(202,50){$\times$}
\put(310,50){$\times$}
\put(0,62){\ldots\ \ $\cdot\cdot\cdot\cdot$ 
111\ \!\!
\ \ \ \ \ \ \ \,$\phi$ \ \ \ \ \ \ \
$\phi$\,\ \ \ \ \ \ \ \ \ \ \ \ \ \ \ \!\ 
$\cdot\cdot\cdot\cdot\cdot\cdot\cdot\,\cdot$
01\ \ \ \ \ \ \ \ \ \ \
\ $\cdot\,\cdot$\ \ \ldots}
\put(38,70){\includegraphics{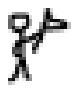}}
\put(36,92){\tiny 1}
\put(261,80){\vector(1,0){20}}
\put(241,75){\includegraphics{queueup.eps}}
\put(238,90){\tiny 2}
\put(251,70){\includegraphics{sticks.eps}}
\put(253,89){\tiny 3}
\put(210,70){\includegraphics{successful.eps}}
\put(207,92){\tiny 4}
\put(80,40){\vector(0,-1){30}}
\put(100,20){Assignment of the ${\bf B}$-symbols to the output
locations}
\put(22,-40){$\times$}
\put(96,-40){$\times$}
\put(135,-40){$\times$}
\put(202,-40){$\times$}
\put(310,-40){$\times$}
\put(22,-30){3255245553\ \ \ \ \ \!\!3451\ \ \ \ 44355521
\ \ \ \ 455251535443551\,\ \ \ \ \ 33\ \ \ldots
}
\put(22,-20){\includegraphics{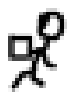}}
\put(29,7){\tiny 1}
\put(40,-6){\vector(1,0){20}}
\put(96,-20){\includegraphics{step4.eps}}
\put(103,7){\tiny 2}
\put(114,-6){\vector(1,0){15}}
\put(135,-20){\includegraphics{step4.eps}}
\put(153,-6){\vector(1,0){20}}
\put(142,7){\tiny 3}
\put(202,-20){\includegraphics{step4.eps}}
\put(220,-6){\vector(1,0){20}}
\put(209,7){\tiny 4}
\put(310,-20){\includegraphics{step4.eps}}
\put(328,-6){\vector(1,0){10}}
\put(317,7){\tiny 5}
\put(45,-117){\framebox(270,66)}
\put(50,-75){\includegraphics{sticks.eps}}
\put(70,-67){{\scriptsize $=$running}}
\put(50,-113){\includegraphics{successful.eps}}
\put(70,-105){{\scriptsize $=$completed simulation}}
\put(180,-75){\includegraphics{queueup.eps}}
\put(200,-67){{\scriptsize $=$queued up}}
\put(180,-113){\includegraphics{step4.eps}}
\put(200,-105){{\scriptsize $=$dispensing output symbols}}
\end{picture}
\vspace{120.0pt}
\caption{Illustration of $\varphi$}
\end{center}
\end{figure}

We now construct the homomorphism $\varphi$ that will be used to prove
Theorem 1. We call the sequence $(x_i)_{i\in\z}$ the {\dof
input sequence}, and the resulting sequence $(\varphi(x)_i)_{i\in\z}$
the {\dof output sequence}. For convenience, we fix a source
probability vector $p=(p(i))_{i\in{\bf A}}$, which we assume has full
support and satisfies $h(p)\ge h(q)+\varepsilon$, and denote ${\bf P}={\bf
   P}_p$. We may also assume without loss of generality that $0\notin
(q(i))_{i\in {\bf B}}$.
For an alphabet ${\bf A}$, denote by
${\bf A}^* = \cup_{n\ge 0} {\bf A}^n$ the set of {\dof finite words
over {\bf A}}. For each $w \in {\bf A}^*$, denote by
$\textrm{length}(w)$ the length of $w$.

The construction is done in several stages. Here is an informal
description of the steps, which are also drawn schematically in Figure 1.

\medskip

{\bf Step 0:} Fix a parameter $t\in\nat$, the {\dof marker length}. Its
value will be some large integer that will be determined later, and will
only depend on the target distribution $q$ and the entropy gap bound
$\varepsilon$.

{\bf Step 1:} Divide the input sequence into {\dof blocks}. A {\dof
marker} is an index $i$ for which
$$ x_i = 2,\ x_{i+1}=x_{i+2}=\ldots=x_{i+t-1}=1. $$ Enumerate the
markers as $\ldots,R_{-2},R_{-1},R_0,R_1,\ldots$, where $R_1$ is the
first marker to the right of the origin. A {\dof block} is the set of
indices between two markers, namely $\{ i : R_{k} < i \le
R_{k+1}\}$. The {\dof input word associated with block $k$} is the
sequence $W_k = (x_i)_{R_{k}+t\le i \le R_{k+1}-1} $, namely the
sequence of input symbols in block $k$, not including the ``211\ldots
1'' patterns.

Under the measure ${\bf P}$, the input words
$\ldots,W_{-1},W_0,W_1,\ldots$ are independent ${\bf A}^*$-valued
random variables. The words $(W_k)_{k\in\z\setminus\{0\}}$ are
identically distributed.  (Note that $W_0$ has a different distribution owing to ``size-biasing'').  All the input words have the property that,
conditioned on the length of $W_k$ being equal to $n$, $W_k$ has
distribution $\tilde{p}_{n,t}$.

{\bf Step 2:} Apply to each input word $W_k$ the function $F_{n,t}$ from
Section 2.2, where $n$ is the length of $W_k$, to obtain a
string $U_k=F_{n,t}(W_k)$ of $N_{n,t}(W_k)$
independent unbiased random bits. $U_k$ is called the {\dof bit string
associated with block $k$}.

{\bf Step 3:} For each block $k\in\z$, attempt to use the random bits in
$U_k$ to simulate a ${\bf B}^{R_{k+1}-R_k}$-valued random variable $B_k$
with distribution $q^{R_{k+1}-R_k}$, using the simulation $(T,S)$ from
Section 2.1. In many blocks, the stopping time will be reached. If the
stopping time is reached, $U_k$ may contain unused bits, which are still
independent and unbiased. For any block $k$ whose stopping time is not
reached, look at $U_{k+1}$ in the next block to the right to find unused
bits to continue the simulation. If now the stopping time is reached,
compute $B_k$. If not, iterate, looking one block further to the right at
each step for unused bits, until the stopping time is reached.  This iteration is done simultaneously for all blocks, in order to maintain translation-equivariance of the construction.  A further complication arises because some bit strings may have length zero, so two or more simulators may try to read the same bits at the same time.  In such situations we give priority to the simulator belonging to the rightmost block.  We will refer to such a situation as a {\dof queue-up}.

 The ergodic
theorem will ensure that, for a proper choice of the parameter
$t$, a.s. enough bits are present overall to complete the process for all
$k\in\z$. Having computed $B_k$, the ${\bf B}$-symbols it contains are
assigned to the indices of the $k$-th block, to produce the output sequence
$\varphi(x)$.

Each block length has exponential tails. For each $k\in \z$, the number of
blocks which must be examined in Step 3 to simulate $B_k$ has exponential
tails. It follows that this homomorphism has exponential tails.

\bigskip
\noindent {\bf Formal definition of $\varphi$}

Let $x=(x_i)_{i\in\z}$. Let $t$, the {\dof marker length}, be a
positive integer to be chosen later.

We first define the marker locations, $(R_k)_{k\in\z}$. Let
$$ R_1 = \min\{ i\ge 0 : x_i=2, x_{i+1}=x_{i+2}=\ldots=x_{i+t-1}
= 1 \}. $$
Inductively, for $k\ge 2$, let
$$ R_k = \min\{ i > R_{k-1} : x_i=2, x_{i+1}=x_{i+2}=\ldots=x_{i+t-1}
= 1\}. $$
Let
$$ R_0 = \max\{ i < 0 : x_i=2, x_{i+1}=x_{i+2}=\ldots=x_{i+t-1}=1 \}. $$
Inductively, for $k\ge 1$, let
$$
R_{-k} = \max\{ i < R_{-k+1} : x_i=2, x_{i+1}=x_{i+2}=\ldots=x_{i+t-1}
=1 \}. $$
It follows from well-known facts of elementary probability theory
that $(R_k)_{k\in\z}$ are defined for ${\bf P}$-almost every input sequence
$(x_i)_{i\in\z}$. For all $k\in \z$, let $\{ i: R_{k} < i \le R_{k+1}
\}$ be the $k$-th {\dof block}. Let $W_k = (x_i)_{R_{k}+t\le i\le
   R_{k+1}-1}$ be the {\dof input word} associated with block $k$, an ${\bf
   A}^*$-valued random variable, and let $L_k =
\textrm{length}(W_k)=R_{k+1}-R_{k}-t$ be its length. Let $\lambda_k =
R_{k+1}-R_k$. Note that by definition, $W_k$ does not contain the pattern
``2 followed by $t-1$ 1's''.

Let $N_{n,t}, F_{n,t}, G_{n,t}$ be as in Section 2.2. For all $k\in
\z$, let $U_k = F_{L_k,t}(W_k)$. $U_k$ is a $\{0,1\}^*$-valued random
variable, called the {\dof bit string} associated with block $k$. Denote
its length by $ V_k = \textrm{length}(U_k) = N_{L_k,t}(W_k)$. Let
$U_k=(\epsilon_k(1),\epsilon_k(2),\ldots,\epsilon_k(V_k))$ be the bits
comprising $U_k$.

For any $\ell\in\nat$, let the pair $(T_\ell,S_\ell)$ be a simulation of
the distribution $q^\ell$ from independent unbiased bits, as in Section
2.1.  (Recall that $T_\ell$ is the stopping time and $S_\ell$ is the output symbol of the simulation).
For an input $x\in\{0,1\}^*$, say that {\bf the simulation
   $(T_\ell,S_\ell)$ is successful for input $x$} if for some
$y\in\{0,1\}^{\nat}$ (and hence also, by the definition of a simulation,
for \emph{all} $y\in\{0,1\}^{\nat}$) we have $T_\ell(x*y)\le
\textrm{length}(x)$, where $x*y$ is the concatenation of $x$ by $y$.

In Step 3, for each $k\in \z$ we will generate a ${\bf B}^*$-valued random
variable $B_k=(\beta_k(1),\ldots,\beta_k(\lambda_k))$ such that,
conditioned on $(\lambda_k)_{k\in\z}$, the $B_k$ are independent and
each $B_k$ has distribution $q^{\lambda_k}$. To do this, in each Step
(3,n) for $n\ge 0$, to each block $k\in\z$ we assign the following: A pair
$(J_k^n,M_k^n)$ with $J_k^n\ge k$ and $1\le M_k^n \le V_{J_k^n}$, called the
{\dof position of the $k$-th simulator at Step $(3,n)$}  (here ``position $(j,k)$'' refers to the $k$th bit of block $j$); a word ${\cal
   Z}_k^n\in \{0,1\}^*$ called the {\dof input read by the $k$-th simulator
   by Step $(3,n)$}; and a set ${\cal G}_k^n$ of pairs $(j, m)\in
\z\times\nat$ such that ${\cal Z}_k^n$ is the concatenation of all the bits
$\epsilon_j(m)$, $(j,m)\in {\cal G}_k^n$, arranged in lexicographical order
on $(j,m)$, called the {\dof set of positions used by the $k$-th simulator
   by Step $(3,n)$}. If for input ${\cal Z}_k^n$, the simulation
$(T_{\lambda_k},S_{\lambda_k})$ is successful, then let $B_k =
S_{\lambda_k}({\cal Z}_k^n*y)$ for some (and hence all)
$y\in\{0,1\}^{\nat}$ and say that $B_k$ {\dof was computed by step
   $(3,n)$}. For a pair $(j,m)$, with $j\in\z$, $V_j>0$ and $1\le m \le
V_j$, denote
$$ \textrm{NEXT}(j,m) = \left\{ \begin{array}{ll}
   (j,m+1) & \textrm{if }m<V_j, \\
   (\min\{j'>j : V_{j'}>0\},1) & \textrm{if }m=V_j. \end{array}\right.
$$
the {\dof next bit position after $(j,m)$} (which is a random
variable).

{\bf Step (3,0):} For all $k\in\z$, set ${\cal G}_k^0 = \emptyset$
(the empty set), ${\cal Z}_k^0 = \phi$ (the empty string). Set $ J_k^0
= \min\{ j \ge k : V_j > 0 \}$ and $M_k^0 = 1$. (It follows easily
from Lemmas \ref{associated} and \ref{tvalue} below that $J_k^0$ are
a.s. defined and finite.) No $B_k$'s are computed.

{\bf Step (3,n):} For each $k\in\z$, if $B_k$ was computed by time
$n-1$, set $(J_k^n,M_k^n)=(J_k^{n-1},M_k^{n-1})$, ${\cal G}_k^n =
{\cal G}_k^{n-1}$ and ${\cal Z}_k^n = {\cal Z}_k^{n-1}$. Otherwise,
check if position $(J_k^{n-1},M_k^{n-1})$ was used by \emph{some}
simulator by Step $(3,n-1)$, i.e. whether $(J_k^{n-1},M_k^{n-1})$ is
in $\cup_{k'\in\z} {\cal G}_{k'}^{n-1}$. If yes: set ${\cal G}_k^n =
{\cal G}_k^{n-1}$, ${\cal Z}_k^n={\cal Z}_k^{n-1}$, and $(J_k^n, M_k^n) =
\textrm{NEXT}(J_k^{n-1},M_k^{n-1})$.

If position $(J_k^{n-1},M_k^{n-1})$ was not used by any simulator by
Step $(3,n-1)$: Check if for some $k'>k$ we have
$(J_k^{n-1},M_k^{n-1})=(J_{k'}^{n-1},M_{k'}^{n-1})$, a phenomenon we
refer to as a {\dof queue-up} of the $k$-th simulator. 
If there is a queue-up, set ${\cal G}_k^n = {\cal G}_k^{n-1}$, ${\cal
   Z}_k^n = {\cal Z}_k^{n-1}$, and $(J_k^n, M_k^n) =
\textrm{NEXT}(J_k^{n-1},M_k^{n-1})$. If the $k$-th simulator is not queued
up: Set ${\cal G}_k^n = {\cal G}_k^{n-1}\cup \{(J_k^{n-1},M_k^{n-1})\}$ and
${\cal Z}_k^n = {\cal Z}_k^{n-1}* \epsilon_{J_k^{n-1}}(M_k^{n-1})$. If now
the simulation $(T_{\lambda_k},S_{\lambda_k})$is successful for input
${\cal Z}_k^n$, set $B_k=S_{\lambda_k}({\cal Z}_k^n * y)$ for some (and
hence all) $y\in\{0,1\}^{\nat}$, set $(J_k^n,M_k^n)=(J_k^{n-1},M_k^{n-1})$,
and say that $B_k$ was computed at Step $(3,n)$. Otherwise, set
$(J_k^n,M_k^n) = \textrm{NEXT}(J_k^{n-1}, M_k^{n-1})$.

We will show later that, if the marker length $t$ is chosen large
enough, then for ${\bf P}$-almost every input sequence
$(x_i)_{i\in\z}$, for all $k\in\z$ there exists an $n\ge 1$ for which
$B_k$ was computed by Step $(3,n)$. So all the $B_k$'s are
a.s. defined ${\bf B}^*$-valued random variables.

Note that it is immediate from the definition that
$\textrm{length}(B_k)=\lambda_k=R_{k+1}-R_k$. Let
$B_k=(\beta_k(1),\beta_k(2),\ldots,\beta_k(\lambda_k))$ be the ${\bf
   B}$-symbols comprising $B_k$. For each $i\in\z$, we define
$(\varphi(x))_i$ as follows. Let $K(i)\in\z$ be the index of the block
containing $i$, namely the unique $k\in\z$ for which $R_{k} < i \le
R_{k+1}$, and set
$$ (\varphi(x))_i = \beta_{K(i)}(i-R_{K(i)}). $$

\bigskip This completes the formal definition of $\varphi$. Figure 1
shows a schematic illustration of the construction. Table 1 summarizes
our main notation for convenient reference.

\newpage

{\linespread{1.2}
\begin{center}
\begin{table}[!h]
\caption{Summary of notation}
\vspace{10.0pt}
\begin{tabular}{l|p{270.0pt}}
Symbol & Meaning \\
\hline
${\bf A}^*$ & finite words over ${\bf A}$ \\
$(T_\ell,S_\ell)$ & simulation of $q^\ell$ from independent unbiased bits \\
$(N_{n,t},F_{n,t}, G_{n,t})$ & simulation of independent unbiased bits
from pattern-free block \\
$\pi$ & ``forbidden pattern'' $(2,1,1,\ldots,1)$\ \  ($t-1$ ones) \\
$E_{n,t}$ & ${\bf A}$-$n$-tuples not containing $\pi$ \\
$x=(x_i)_{i\in\z}$ & input sequence \\
${\bf A}=\{j : 1 \le j \le a\}$ & source alphabet \\
${\bf B}=\{j : 1 \le j \le b\}$ & target alphabet \\
$p=(p(j))_{1\le i\le a}$ & source distribution \\
$q=(q(j))_{1\le j\le b}$ & target distribution \\
$\varepsilon$ & lower bound on the entropy gap $h(p)-h(q)$ \\
$\tilde{p}_{n,t}$ & $p^n$ conditioned on $E_{n,t}$ \\
$t$ & marker length \\
$R_k$ & $k$-th marker \\
$\lambda_k=R_{k+1}-R_k$ & length of $k$-th block \\
$W_k=(x_i)_{R_k+t \le i \le R_{k+1}-1}$ & $k$-th input word  \\
$L_k=\lambda_k-t$ & length of $W_k$ \\
$U_k=(\epsilon_k(1),\ldots,\epsilon_k(V_k))$ 
  & bit string associated with block $k$ \\
$V_k$ & length of $U_k$ \\
$(J_k^n,M_k^n)$ & position of the $k$-th simulator at Step $(3,n)$ \\
${\cal Z}_k^n$ & input read by the $k$-th simulator by Step $(3,n)$ \\
${\cal G}_k^n$ & positions used by the $k$-th simulator by Step
$(3,n)$ \\
$\textrm{NEXT}(j,m)$ & next bit position after $(j,m)$ \\
$B_k=(\beta_k(1),\ldots,\beta_k(\lambda_k))$ & ${\bf B}^*$-valued
r.v. computed by the $k$-th simulator \\
$K(i)$ & index of block containing $i$ \\
$\varphi(x)=(\varphi(x)_i)_{i\in\z}$ & output sequence
\end{tabular}
\end{table}
\end{center}
}

\section{Proof of Theorem 1}

%

\begin{lemma} \label{wordsindep}
The ${\bf A}^*$-valued random variables $(W_k)_{k\in\z}$
are independent. For each $k\in\z$ and $n\in\{0,1,2,\ldots\}$, $W_k$
conditioned on $\{\textrm{length}(W_k)=n\}$ has distribution
$\tilde{p}_{n,t}$. The non-central block lengths
$(R_{k+1}-R_k)_{k\in\z\setminus\{0\}}$ are identically distributed
with exponential tails, and the central block length \hbox{$R_1-R_0$} has
exponential tails.
\end{lemma}

\begin{pf} Let $\mu_0, \mu_1$ be the measures on ${\bf A}^*$ defined as
follows.
\begin{eqnarray*}
\lefteqn{\mu_0\big(\{(a_1,a_2,\ldots,a_n)\}\big)} \\
&=& \left\{ 
\begin{array}{ll}
(n+t) p(2)^2 p(1)^{2(t-1)} \prod_{j=1}^n p(a_i), & (a_1,\ldots,a_n) \in
E_{n,t}\\ 0 & \textrm{otherwise} \end{array}\right.
\\ 
\lefteqn{\mu_1\big(\{(a_1,a_2,\ldots,a_n)\}\big)} \\
&=& \left\{ 
\begin{array}{ll}
p(2)p(1)^{t-1} \prod_{j=1}^n p(a_i),\qquad\qquad\ 
 & (a_1,\ldots,a_n) \in
E_{n,t}\\ 0 & \textrm{otherwise} \end{array}\right.
\end{eqnarray*}
We claim that for any $j\ge 0$ and $(w_k)_{-j\le k\le j} \subset {\bf
A}^*$,
\begin{equation}\label{eq:indep}
{\bf P}\bigg((W_k)_{-j\le k\le j}=(w_k)_{-j\le k\le j}\bigg) = \mu_0(w_0)
\prod_{-j\le k\le j,\  k\ne 0} \mu_1(w_k).
\end{equation}
This will prove that $(W_k)_{k\in\z}$ are independent, with $W_0$
having distribution $\mu_0$ and all the other $W_k$'s having
distribution $\mu_1$ -- with both $\mu_i$'s clearly having the desired
property that conditioning on the length $n$ gives
$\tilde{p}_{n,t}$. (Incidentally, it will also prove that $\mu_0,
\mu_1$ are probability measures, although this can be checked
directly.) Indeed, to prove \eqref{eq:indep}, let $\pi=211\ldots 1$ be
the word ``2 followed by $t-1$ 1's'', and let $w$ be the word obtained
by the concatenation
$$ w = \pi * w_{-j} * \pi * w_{-j+1} * \ldots * \pi * w_j * \pi. $$
Denote $\textrm{len}^+ = \sum_{1 \le k \le j} \textrm{length}(w_k)$,
$\textrm{len}^- = \sum_{-j \le k \le 0} \textrm{length}(w_k)$.
Then, if all the $w_k$'s do not contain the pattern $\pi$, we have
\begin{multline}\label{eq:eventeq}
\bigg\{ (W_k)_{-j\le k\le j}=(w_k)_{-j\le k\le j} \bigg\} \\
= \bigcup_{r=0}^{\textrm{length}(w_0)+t-1}
\bigg\{ (x_i)_{i=r - (j+1)t-\textrm{len}^-}
^{r + (j+1)t+ \textrm{len}^+} = w \bigg\},
\qquad
\end{multline}
and furthermore the above union is disjoint.
(In words, this simply means that $(W_k)_{-j\le k\le j}=(w_k)_{-j\le
k\le j}$ if and only if for some $r$, a string of $x_i$'s centered
around the origin such that the offset of the $0$-th block relative to
the origin is equal to $r$, is equal to $w$. This follows directly
from the definitions.) Therefore
\begin{multline*}
{\bf P}\bigg((W_k)_{-j\le k\le j}=(w_k)_{-j\le k\le j}\bigg) =
\sum_{r=0}^{\textrm{length}(w_0)+t-1} {\bf P}\bigg(
(x_i)_{i=r - (j+1)t-\textrm{len}^-}
^{r + (j+1)t+ \textrm{len}^+} = w \bigg) \\
= (\textrm{length}(w_0)+t)\cdot p^{\textrm{length}(w)}(w) =
\mu_0(w_0) \prod_{-j\le k\le j,\  k\ne 0} \mu_1(w_k).
\end{multline*}
If some $w_k$ contains $\pi$, the event on the left hand side of
\eqref{eq:eventeq} is empty, and \eqref{eq:indep} holds trivially.

It remains to prove that the block lengths $R_{k+1}-R_k$ have
exponential tails, or equivalently to prove the same for
$L_k=\textrm{length}(W_k) = R_{k+1}-R_k-t$. Denote $c =
p^t(\pi)=p(2)p(1)^{t-1}$. Then for $k\neq 0$,
\begin{multline*}
{\bf P}\bigg(L_k=n\bigg) = \sum_{w\in E_{n,t}}
\mu_1(w) = c \sum_{w \in E_{n,t}} p^n(w) \\
= c\,{\bf P}\bigg( (x_i)_{i=0}^{n-1}\textrm{ does not contain }\pi
\bigg) \qquad\qquad\qquad\qquad\qquad
\\
\le c\,{\bf P}\bigg( (x_i)_{i=j t}^{jt+t-1} \neq \pi,\ \
j=1,2,\ldots,\left\lfloor \frac{n-t}{t}\right\rfloor \bigg)=
c (1-c)^{\lfloor (n-t)/t\rfloor},
\end{multline*}
which decays exponentially in $n$. Similarly, since the
distribution $\mu_0$ is at most a factor $O(n)$ times $\mu_1$,
$L_0=\textrm{length}(W_0)$ also has exponential tails.
\end{pf}

\begin{lemma} \label{associated}
The associated bit strings $(U_k)_{k\in\z}$ are independent
$\{0,1\}^*$-valued random variables. For each $k\in\z$ and
$n\in\{0,1,2,\ldots\}$, $U_k$ conditioned on $\{V_k = n\}$ has the
uniform distribution on $\{0,1\}^n$. $(U_k)_{k\in\z\setminus\{0\}}$
are identically distributed. $(V_k)_{k\in\z}$ have exponential tails.
\end{lemma}

\begin{pf} It is immediate from Lemma \ref{wordsindep} that $(W_k)_{k\in\z}$ are
independent and $(W_k)_{k\in\z\setminus\{0\}}$ are identically
distributed, therefore the same holds for the sequence $(U_k)$, as required. Now, for any $k\in\z$, $V_k = N_{L_k,t}(W_k) \le (\log a/\log
2)L_k$ by Theorem \ref{simu2}(iii), so it has exponential tails by
Lemma \ref{wordsindep}. Finally, let $k\in\z$ and
$n\in\{0,1,2,\ldots\}$.  Then by Theorem \ref{simu2} and Lemma
\ref{wordsindep}, for any $w\in\{0,1\}^n$,
\begin{multline*}
{\bf P}\bigg(U_k = w\ \bigg|\ V_k = n\bigg) \\ =
\sum_m {\bf P}\bigg(L_k=m\ \bigg|\ V_k=n\bigg)
       {\bf P}\bigg(U_k = w\ \bigg|\ V_k = n,\ L_k=m\bigg) 
\qquad\qquad\qquad\ \ \ \\
= \sum_m {\bf P}\bigg(L_k=m\ \bigg|\ V_k=n\bigg) \qquad\qquad\qquad\qquad\qquad\qquad\qquad\qquad\qquad\qquad\qquad\\
\qquad\qquad\qquad\times{\bf P}\bigg(F_{m,t}(W_k)= w\ \bigg|\ \textrm{length}(W_k)=m,\
N_{m,t}(W_k)=n \bigg) \qquad\qquad\qquad\qquad\qquad\qquad\qquad\qquad\qquad\qquad\qquad\\
= \sum_m {\bf P}\bigg(L_k=m\ \bigg|\ V_k=n\bigg)\cdot 2^{-n} = 2^{-n}.
\qquad\qquad\qquad\qquad\qquad\qquad\qquad\qquad
\end{multline*}
\end{pf}

\begin{lemma} \label{tvalue}
The marker length $t$ may be chosen so that for $k\ne 0$, 
\begin{equation}\label{eq:entropy}
{\bf E}(V_k) > {\bf E}(T_{\lambda_1}).
\end{equation}
\end{lemma}

\begin{pf} 
Consider temporarily a new probability measure $\tilde{\bf P}$, with expectation operator $\tilde{\bf E}$, under which 
$W_0$ 
has the same distribution as
$W_1$ and is independent of all other random variables, while all other random variables have the same distribution as before.  By Lemma
\ref{wordsindep} and the computation in its proof, the
process $(W_k)_{k\in\z}$ is now isomorphic to the induced dynamical system
$B(p)_{\big| {\cal M}}$, where
$${\cal M}=\{ x\in {\bf A}^\z : (x_i)_{i=0}^{t-1} = \pi \}. $$
(Recall that $\pi=(2,1,...,1)$ is the forbidden pattern.)
By Abramov's
formula (\cite{petersen}, p. 257--259), this dynamical system has entropy
$h(W_1)=h(p)/p^t(\pi)$. This is also equal to $h(p)\tilde{\bf E}(\lambda_1)=h(p){\bf E}(\lambda_1)$,
since by Kac's formula (\cite{petersen}, p. 46), the expected return time
(or expected block length) is the reciprocal of the probability of the
inducing set.

By Theorem \ref{simu2}(ii), the mapping $W_1 \mapsto (U_1, L_1,
G_{L_1,t}(W_1))$ is injective. Therefore, using Lemma 5 and elementary
properties of the entropy function,
\begin{eqnarray*}
h(p){\bf E}(\lambda_1) &=& h(W_1)=h(U_1, L_1, G_{L_1,t}(W_1)) \le h(U_1) +
h(L_1, G_{L_1,t}(W_1)) \\ &=& h(U_1 | V_1) + h(V_1) + h(L_1, G_{L_1,t}(W_1))
\\ &=& {\bf E}(V_1)\log 2 + h(V_1) + h(L_1, G_{L_1,t}(W_1)).
\end{eqnarray*}
By Theorem \ref{simu1}, we have
$$ {\bf E}(T_{\lambda_1})\log 2 \le 6\log 2 + {\bf E}(\lambda_1)h(q) \le
6\log 2 + {\bf E}(\lambda_1)(h(p)-\varepsilon). $$
Combining the last two inequalities gives that
$$
{\bf E}(V_1) - {\bf E}(T_{\lambda_1}) \ge \frac{\varepsilon}{\log 2}
{\bf E}(\lambda_1) - 6 - \frac{h(V_1) + h(L_1,G_{L_1,t}(W_1))}{\log 2}. $$
To bound the negative terms on the right-hand side, recall the following
properties of the entropy of integer-valued random variables (see
\cite{coverthomas}, Lemma 12.10.2): If $X$ is a random variable with finite
expectation that takes values in $\nat$, then $h(X) \le h(Y)$, where $Y$
has a geometric distribution with ${\bf E}(Y) = {\bf E}(X)$. Furthermore,
$h(Y)=O(\log {\bf E}(Y))$ when the expectation is large. This implies,
using the fact that $V_1 \le (\log a/\log 2)L_1$, that for some positive
constant $C$ (that depends on $a$),
\begin{eqnarray*} 
h(L_1) &\le& C \log {\bf E}(\lambda_1), \\
h(V_1) &\le& C \log {\bf E}(\lambda_1), \\
h(L_1, G_{L_1,t}(W_1)) &=& h(L_1) + h(G_{L_1,t}(W_1)\ |\  L_1) 
\\ &\le& C\log {\bf
   E}(\lambda_1) + (a-1)\log {\bf E}(\lambda_1),
\end{eqnarray*}
(because the range of $G_{n,t}$ is
$\{1,2,\ldots,n^{a-1}\}$). Therefore
$${\bf E}(V_1) - {\bf E}(T_{\lambda_1}) \ge f({\bf E}(\lambda_1)) $$
for some function
$$ f(u) = f_{\varepsilon}(u) = \frac{\varepsilon}{\log 2}u - 6 - C' \log u, $$
where $C'>0$ is a constant that depends only on $a$. Now, $f(u)\to\infty$ as
$u\to\infty$. Choosing the marker length $t$ sufficiently high will force
${\bf E}(\lambda_1)=1/p^t(\pi)$ to be large, uniformly over all probability
vectors $p$ under consideration, i.e. that satisfy 
$h(p)\ge h(q)+\varepsilon$ and $p(i)>0$ for all $i \in {\bf A}$
 (note that $p(1)$ is bounded away from 1
because the alphabet size is fixed and $h(p)$ is bounded away from 0). So
for sufficiently large $t$ we get ${\bf E}(V_1)-{\bf E}(T_{\lambda_1})>0$,
proving the lemma.
\end{pf}

From now on we consider $t$ as having a fixed value for which
\eqref{eq:entropy} holds. Note that in particular it follows from
Lemma \ref{tvalue} that almost surely, $V_k>0$ for infinitely many
positive values of $k$. Therefore $(J_k^n,M_k^n)$ in Steps $(3,0)$,
$(3,n)$ are a.s. defined.

%
%

Recall the notion of {\dof stochastic domination}. Let $\Lambda$ be
a compact metric space on which there is defined a partial order
``$\preceq$'', and assume that $\preceq$ is a closed subset of $\Lambda
\times \Lambda$.
If $X,Y$ are two
random variables (not necessarily defined on the same probability
space) taking values in $\Lambda$,
denote $X \preceq_{\textrm{stoc}} Y$ (read: {\dof ``$X$ is
stochastically dominated by $Y$''}) if for any $e\in\Lambda$ we have
$$ {\bf P}(X \succeq e) \le {\bf P}(Y \succeq e). $$
It is known (see \cite[Th. 2.4, p. 71]{liggett})
that, under the above topological assumptions,
$X \preceq_{\textrm{stoc}} Y$ if and
only if there exist random variables $X',Y'$, defined \emph{on the
same probability space}, such that the variable
$X'$ has the same distribution as
$X$, the variable $Y'$ has the same distribution as $Y$, and
$$ {\bf P}(X' \preceq Y') = 1. $$

On the set $\{0,1\}^{\#}=
\{0,1\}^* \cup \{0,1\}^{\mathbb{N}}$ of all finite and
infinite bit strings, let $w \preceq w'$
denote the order ``$w$ is a prefix of $w'$''. Equip $\{0,1\}^{\#}$
with the topology consisting of open sets of the form
$\left\{ w' \in \{0,1\}^{\#} : w \preceq w' \right\}$ for
$w \in \{0,1\}^*$. It is not difficult to verify that $\{0,1\}^{\#}$
with this topology is a compact metric space, and that $\preceq$ is
a closed subset of $\{0,1\}^{\#}\times \{0,1\}^{\#}$.

On the set
$(\{0,1\}^{\#})^\z$ of $\z$-indexed vectors of 
finite and infinite bit strings, let
$\preceq^\z$ denote the (strong) product order
of $\preceq$, i.e., we define
$$ (w_k)_{k\in\z} \preceq^\z (w_k')_{k\in\z}\ \ \iff\ \ 
w_k \preceq w_k'\ \ \textrm{for all }k\in\z. $$

For each $\ell\in\nat$, let $Z_\ell$ be the following $\{0,1\}^*$-valued
random variable: take independent unbiased bits $\nu_1,\nu_2,\nu_3,\ldots$,
and set
$$
Z_\ell = (\nu_1,\nu_2,\ldots,\nu_{T_\ell(\nu_1,\nu_2,\ldots)}), $$
where
$T_\ell$ is the stopping time of the simulation of $q^\ell$ from
independent unbiased bits as in Section 3. We call $Z_\ell$ (an instance of)
{\dof an acceptable input for the simulation $(T_\ell,S_\ell)$}. Let
$(Z_{\ell,k})_{\ell\in\nat, k\in\z}$ be an infinite array of independent
random variables, where $Z_{\ell,k}$ has the same distribution as
$Z_\ell$. Denote $\tau_{\ell,k}=\textrm{length}(Z_{\ell,k})$. Clearly
$\tau_{\ell,k}$ is equal in law to $T_\ell(\nu_1,\nu_2,\ldots)$.

\begin{lemma} \label{induction}
For each $n=0,1,2,\ldots$ we have
\begin{mylist}
\item ${\cal Z}_k^n$ is the concatenation of the bits $\epsilon_j(m)$
for all $(j,m)\in{\cal G}_k^n$, arranged by lexicographical order on
$(j,m)$.
\item $({\cal G}_k^n)_{k\in\z}$ are disjoint sets.
\item The conditional distribution of $({\cal Z}_k^n)_{k\in\z}$ given
$(R_k)_{k\in\z}=(r_k)_{k\in\z}$ is a.s. stochastically dominated, in the order
$\preceq^\z$, by $(Z_{r_{k+1}-r_k,k})_{k\in\z}$.
\item The conditional distribution of $(\textrm{length}({\cal
Z}_k^n))_{k\in\z}$ given $(R_k)_{k\in\z}=(r_k)_{k\in\z}$ is a.s. stochastically
dominated, in the order $\le^\z$ (the product order of the usual order
on numbers), by $(\tau_{r_{k+1}-r_k,k})_{k\in\z}$.
\end{mylist}
\end{lemma}

\begin{pf} Claim (i) follows trivially by induction on $n$, as the
simulator locations $(J_k^n,M_k^n)$ are obviously increasing in the
lexicographical order.

Claim (ii) follows by induction on $n$, by noting that ${\cal G}_k^n$
is always obtained from ${\cal G}_k^{n-1}$ by the addition of at most
one location $(j,m)$, and, since we made allowance for the phenomenon
of queue-ups, where multiple simulators are at the same location
during Step $(3,n)$, any given location $(j,m)$ is added to ${\cal
G}_k^n$ for at most one value of $k\in\z$, keeping the ${\cal
G}_k^n$'s disjoint.

It remains to prove Claim (iii), which implies (iv) trivially. To do this,
let $(\theta_{k,j})_{k\in\z, j\in\nat}$ be an array of independent unbiased
bits which are independent of all other random variables. For each $k\in\z$
define ${\cal X}_k = {\cal Z}_k^n * (\theta_{k,j})_{j\ge 1}$. Because of
Lemma \ref{associated} together with claims (i) and (ii) proven above, it
follows that, conditional on $(R_k)_{k\in\z}=(r_k)_{k\in\z}$,
we have that
$({\cal X}_k)_{k\in\z}$ is a sequence of independent infinite sequences of
independent unbiased bits. Set $\xi_k = ({\cal X}_{k,j})_{1\le j \le
   T_{r_{k+1}-r_k}({\cal X}_k)}$. Then (still working conditionally) 
$(\xi_k)_{k\in\z}$ are independent and for each $k\in\z$, $\xi_k$ has the
distribution of $Z_{r_{k+1}-r_k}$. Also, from the construction necessarily
${\cal Z}_k^n \preceq \xi_k$. We have constructed a realization of 
$(Z_{r_{k+1}-r_k,k})_{k\in\z}$ that dominates $({\cal Z}_k^n)_{k\in\z}$,
thereby proving the stochastic domination claim.
\end{pf}

We shall use the following {\dof mass-transport lemma}:

\begin{lemma}\label{masstrans}
Let $f:\z\times\z\to\mathbb{R}$ satisfy $f(x+c,y+c)=f(x,y)$ for all
$x,y,c\in\z$. Then for all $x\in\z$,
$$ \sum_{y\in\z} f(x,y) = \sum_{y\in\z} f(y,x). $$
\end{lemma}

\begin{pf} Taking $c=x-y$ gives $f(x,y)=f(2x-y,x)$, so
$$ \sum_{y\in\z} f(x,y) = \sum_{y\in\z} f(2x-y,x) = \sum_{u\in\z}
f(u,x). $$
\end{pf}

\begin{lemma}\label{bkdefined}
Almost surely, for all $k\in\z$ there exists an $n\ge 1$ for which
$B_k$ was computed by Step $(3,n)$.
\end{lemma}

\begin{pf} As in the proof of Lemma \ref{tvalue}, we introduce the probability measure $\tilde{\bf P}$ with expectation operator $\tilde{\bf E}$ under which $(W_k)_{k\in{\mathbb Z}}$ are i.i.d.  Since ${\bf P}$ is absolutely continuous with respect to
    $\tilde{\bf P}$ (see the proof of Lemma \ref{wordsindep}), it is enough
    to prove that each $B_k$ is eventually computed, $\tilde{\bf P}$-a.s.

For each $k\in\z$, define ${\cal Z}_k^\infty = \lim_{n\to\infty} {\cal
   Z}_k^n$ (a possibly infinite bit sequence) and ${\cal G}_k^\infty =
\lim_{n\to\infty} {\cal G}_k^n = \cup_{n=1}^\infty {\cal G}_k^n$. Define
$f:\z\times\z \to \mathbb{R}$ by
\begin{eqnarray*} f(k,j) &=& {\tilde{\bf E}}\bigg(
\big| \big\{ (j,i) : i\in\z,\ \ 1 \le i \le V_j \big\}
 \cap {\cal G}_k^\infty\big| \bigg)
\\ &=&
\tilde{\bf E}\big(\textrm{number of bits from $U_j$ read by the
   $k$-th simulator}\big).
\end{eqnarray*}
The function $f$ satisfies the assumption of Lemma \ref{masstrans}, since
$(W_k)_{k\in\z}$ is a stationary sequence under $\tilde{\bf P}$, and
it is easy to see from the construction that $(U_k,V_k,{\cal Z}_k^\infty,
{\cal G}_k^\infty)_{k\in\z}$ are all generated from $(W_k)_{k\in\z}$ in a
shift-equivariant manner. Therefore for all $k\in\z$, we have
\begin{equation}\label{eq:transport}
\sum_{j\in\z} f(k,j) = \sum_{j\in\z} f(j,k).
\end{equation}
Denote the quantity in \eqref{eq:transport} by $g$ (by stationarity it does
not depend on $k$). The left-hand side of \eqref{eq:transport} is equal to
$\tilde{\bf E}(|{\cal G}_k^\infty|)$, the
expected number of bits eventually read by the $k$-th simulator. By Lemma
\ref{induction}(iv), $g \le {\bf E}(T_{\lambda_1})$. The
right-hand side is equal to the total number of bits from the $k$-th
associated bit string $U_k$ eventually used by \emph{any} simulator, and
clearly cannot exceed $\tilde{\bf E}(V_k) = {\bf E}(V_1)$. In particular, this implies that $g<\infty$, so almost surely
${\cal Z}_k^\infty$ is a finite string.

Assume that with positive $\tilde{\bf P}$-probability, $B_k$ is not
computed by Step $(3,n)$ for any $n$. Since ${\cal Z}_k^\infty$ is finite,
the only way for this to happen is for the $k$-th simulator to eventually
fail to find any unused bits in the blocks to its right. Because of
stationarity, by the ergodic theorem
this implies that a.s., this happened to a positive
proportion of the simulators. Therefore, a.s. \emph{all} the bits are
eventually used! In other words, there is the equality $g={\bf E}(V_1)$. We
have shown ${\bf E}(V_1) = g \le {\bf E}(T_{\lambda_1})$, in contradiction
to Lemma \ref{tvalue}. So the assumption that $B_k$ was not computed with
positive probability is false.
\end{pf}

\begin{lemma}\label{bkdist}
Conditioned on $(R_k)_{k\in\z}=(r_k)_{k\in\z}$, the random variable
$(B_k)_{k\in\z}$ are independent, and for each $k\in\z$,
$B_k$ has distribution $q^{r_{k+1}-r_k}$.
\end{lemma}

\begin{pf} By Lemma \ref{induction}, conditioned on
    $(R_k)_{k\in\z}=(r_k)_{k\in\z}$, the random vector
$({\cal Z}_k^\infty)_{k\in\z}$ is
    stochastically dominated in the product prefix order by
    $(Z_{r_{k+1}-r_k,k})_{k\in\z}$. Assume that these two sequences are
    defined on the same space and that there is actual a.s. domination. For
    each $k\in\z$, both ${\cal Z}_k^\infty$ and $Z_{r_{k+1}-r_k,k}$ are
    acceptable inputs for the simulation
    $(T_{r_{k+1}-r_k},S_{r_{k+1}-r_k})$. So, since
${\cal Z}_k^\infty \preceq Z_{r_{k+1}-r_k,k}$, necessarily
${\cal Z}_k^\infty = Z_{r_{k+1}-r_k,k}$. We have shown that
$({\cal Z}_k^\infty)_{k\in\z} = (Z_{r_{k+1}-r_k,k})_{k\in\z}$. Therefore
(still working conditionally) $({\cal Z}_k^\infty)_{k\in\z}$ are
independent and each ${\cal Z}_k^\infty$ has the distribution of an
acceptable input for the simulation
$(T_{r_{k+1}-r_k},S_{r_{k+1}-r_k})$. Therefore, since
$B_k=S_{r_{k+1}-r_k}({\cal Z}_k^\infty)$,
$(B_k)_{k\in\z}$ are independent, and for each $k\in\z$,
$B_k$ has distribution $q^{r_{k+1}-r_k}$.
\end{pf}

\begin{lemma} \label{phidist}
$((\varphi(x))_i)_{i\in\z}$ are i.i.d. with distribution $q$.
\end{lemma}

\begin{pf}
The mapping $i\to (K(i),i-R_{K(i)})$ is obviously
injective. This means that the $(\varphi(x))_i=\beta_{K(i)}(i-R_{K(i)})$ get
assigned different $\beta_k(j)$ symbols. Conditioned on $(R_k)_{k\in\z}$,
these symbols are all independent ${\bf B}$-symbols with distribution $q$.
This immediately implies the claim of the lemma.
\end{pf}

\begin{lemma} \label{equivariant}
$\varphi$ is translation-equivariant.
\end{lemma}

\begin{pf} Let $x\in{\bf A}^\z$. Denote $x'={\cal T}(x)$. Our goal is
to prove that for all $i\in\z$, $(\varphi(x'))_i =
(\varphi(x))_{i+1}$.

If $X$ is any of the various quantities in Table 1 which are
implicitly dependent on $x$, denote by $X'$ the corresponding quantity
taken as a function of $x'$ rather than $x$. Consider separately two
cases:

{\bf Case 1:} $R_1 > 0$. In this case, it is easy to check directly
that for all $k\in\z$,
$$ R_k' = R_k-1,\ \  W_k' = W_k,\ \  L_k'=L_k,\ \  U_k'=U_k,\ \  
V_k'=V_k$$
(the markers are shifted by one).
By induction on $n$, for all $k\in\z$ and all $n\ge 0$,
$$ ((J_k^n)',(M_k^n)') = (J_k^n, M_k^n),\ \ ({\cal Z}_k^n)'={\cal Z}_k^n,\ \
({\cal G}_k^n)' = {\cal G}_k^n. $$
Therefore, for all $k\in\z$, $B_k'=B_k$. Furthermore, $K(i)' = K(i+1)$.
Therefore 
$$ (\varphi(x'))_i = \beta_{K(i)'}'(i-R_{K(i)'}') =
\beta_{K(i+1)}(i-(R_{K(i)+1}-1))
 = (\varphi(x))_{i+1}. $$

{\bf Case 2:} $R_1=0$. In this case, check that for all $k\in\z$,
$$ R_k' = R_{k+1}-1,\ \  W_k' = W_{k+1},\ \  L_k'=L_{k+1},\ \  
U_k'=U_{k+1},\ \ V_k'=V_{k+1}$$
(the markers are shifted by 1, and the indexing of the blocks is
shifted by 1). Therefore, by induction
on $n$, for all $k\in\z$ and for all $n\ge 0$,
$$ ((J_k^n)',(M_k^n)') = (J_{k+1}^n, M_{k+1}^n),\ \ 
({\cal Z}_k^n)'={\cal Z}_{k+1}^n,\ \ ({\cal G}_k^n)' = {\cal G}_{k+1}^n. $$
Therefore, for all $k\in\z$, $B_k'=B_{k+1}$. Check as before that now
$ K(i)' = K(i+1)-1$.
Therefore
\begin{eqnarray*}
(\varphi(x'))_i &=& \beta_{K(i)'}'(i-R_{K(i)'}') \\&=&
\beta_{(K(i+1)-1)+1}(i-(R_{(K(i+1)-1)+1}-1))
= (\varphi(x))_{i+1}. \end{eqnarray*}
\end{pf}

The following facts concerning random variables with exponential tails
will be useful.

\begin{lemma} \label{expotails} \ 
\begin{mylist}
\item If $X,Y$ are real-valued random variables with exponential
tails, then $X+Y$ has exponential tails.
\item If $X_1, X_2, \ldots$ are random variables with uniformly
exponential tails (i.e., $(\sup_k {\bf P}(|X_k|\ge n))_{n\in{\mathbb
N}}$ decays exponentially), and $T$ is an $\nat$-valued random variable
with exponential tails, then the random variable $X_T$ has exponential
tails.
\item If $X_1, X_2, \ldots$ are i.i.d. random variables with
exponential tails and mean $\mu$, then, for any $c > 0$, the sequence
$$\bigg(\py\bigg(\bigg|\sum_{i=1}^n X_i-n\mu\bigg|
\geq nc\bigg)\bigg)_{n \in {\mathbb Z}_+}$$
decays exponentially.
\item Suppose $X_1, X_2, \ldots$ are i.i.d. random variables with
exponential tails, and $T$ is an $\nat$-valued random variable with
exponential tails. Denote $S_m = \sum_{k=1}^m X_k$. Then $S_T :=
\sum_{k=1}^T X_k$ has exponential tails.
\end{mylist}
\end{lemma}

\begin{pf} Proof of (i): we have
$ {\bf P}(|X+Y|\ge n) \le {\bf P}(|X|\ge n/2) + {\bf P}(|Y|\ge n/2)$, which decays
exponentially.

Proof of (ii): we have
\begin{eqnarray*}
{\bf P}(|X_T|\ge n) &\le& {\bf P}(T\ge n) + {\bf P}\left(
    \bigcup_{k=1}^n \{ |X_k|\ge n \} \right)
\\ &\le& {\bf P}(T\ge n) + n \sup_{1\le k\le n}{\bf P}(|X_k|\ge n),
\end{eqnarray*}
which decays exponentially.

Claim (iii) is a standard fact from large deviation theory; see for example (\cite{kallenberg} Corollary 27.4).

Proof of (iv): Let $c>0$ and $0<d<1$ be such that ${\bf P}(T\ge n)\le
c\cdot d^n$ for all $n$. Denote $a=1/(2{\bf E}|X_1|)$. Then
\begin{eqnarray*}
{\bf P}(|S_T|\ge n) &=& \sum_{m=1}^\infty 
{\bf P}\big(T=m,\ |S_m|\ge n\big) \le
\sum_{m=1}^\infty {\bf P}\bigg(T=m, \sum_{k=1}^m |X_k| \ge n \bigg) \\
&\le & \sum_{m=1}^{\lfloor a\cdot n \rfloor}
{\bf P}\bigg(T=m, \sum_{k=1}^{\lfloor a\cdot n\rfloor} |X_k| \ge n
\bigg) + \sum_{m={\lfloor a\cdot n\rfloor+1}}^\infty
{\bf P}\bigg(T=m\bigg) \\ &\le&
e n\cdot {\bf P}\bigg(\bigg| \sum_{k=1}^{\lfloor a\cdot n\rfloor}
|X_k| - \lfloor a\cdot n\rfloor {\bf E}|X_1| \bigg| > \frac{n}{2}
\bigg) + 
{\bf P}\bigg( T > a n \bigg)
\end{eqnarray*}
In the last bound, the second term decays exponentially in
$n$. The first term decays exponentially, by (iii) above.
\end{pf}

For each $k\in\z$, let $J_k^\infty = \lim_{n\to\infty} J_k^n$ be the
value of $J_k^n$ for that $n$ for which $B_k$ was computed at Step
$(3,n)$; that is, the index of the rightmost block used by simulator $k$.

\begin{lemma} \label{jkinfty}
$J_0^\infty$ has exponential tails.
\end{lemma}

\begin{pf}
Let $\tau$ be that $n>0$ for which $B_0$ was computed at
Step $(3,n)$. If we define $I_0=0$,
$$ I_1 = \min\{ j \ge 0 : V_j > 0 \}, $$
and inductively for $k>1$,
$$
I_k = \min\{ j > I_{k-1} : V_j > 0 \}, $$
then because of Lemma
\ref{associated}, it is easy to prove in a manner similar to the proof of
Lemma \ref{wordsindep} that $(I_{k+1}-I_k)_{k\ge 0}$ are independent,
$(I_{k+1}-I_k)_{k\ge 1}$ are identically distributed and have the geometric
distribution $Geom(p_0)$, where $p_0={\bf P}(V_1>0)$ (so in particular have
exponential tails), and $I_1$ is stochastically dominated by $(I_2-I_1)+1$.

We know that
$$ J_0^\infty \le I_{\tau} = \sum_{j=1}^\tau (I_j-I_{j-1}). $$
This is because $J_0^0=I_1$, and at any Step $(3,n)$,
$1\le n < \tau$, we have
$(J_0^n,M_0^n)=\textrm{NEXT}(J_0^{n-1},M_0^{n-1})$, so by induction,
$J_0^n \le I_{n+1}$. In particular
$J_0^\infty=J_0^\tau=J_0^{\tau-1}\le I_\tau$.

It follows, by Lemma \ref{expotails}(iv), that it is enough to prove
that $\tau$ has exponential tails. Let $c>0, \delta>0$ be parameters
whose value we will fix shortly. Write
\begin{eqnarray*}
{\bf P}\bigg( \tau > n \bigg) &=& {\bf P}\bigg( B_0\textrm{ undefined by
Step }(3,n)\bigg)
\qquad\qquad\qquad\qquad
\end{eqnarray*}

\vspace{-20.0pt}
\begin{eqnarray}\label{eq:twoterms}
&=& {\bf P}\bigg( \tau > n,\ \ \textrm{length}({\cal Z}_0^n) < c n
\bigg) + {\bf P}\bigg( \tau > n,\ \ \textrm{length}({\cal Z}_0^n) \ge c
n \bigg).\qquad
\end{eqnarray}
By Lemma \ref{induction}(iv), the second term is at most ${\bf
   P}(T_{\lambda_1} \ge cn)$. Note that $T_{\lambda_1}$ has exponential
tails, since by \ref{simu1}(i), ${\bf P}(T_k\ge n)\le \frac{b^k+1}{2^n}$,
whence
$$ {\bf P}(T_{\lambda_1}\ge n) \le {\bf P}(\lambda_1 \ge dn) +
{\bf P}\left(\bigcup_{k\le dn} \{ T_k\ge n \} \right) \le 
{\bf P}(\lambda_1 \ge dn) + \frac{b^{dn}+1}{2^n} $$
can be seen to decay exponentially by taking $d = \log2/(2\log b)$.
It remains to deal with the first term in \eqref{eq:twoterms},
$$ {\bf P}(E_n) :=
{\bf P}\bigg( \tau > n,\ \ \textrm{length}({\cal Z}_0^n) < c n
\bigg).
$$
Note the event inclusion
$$ E_n \subseteq \bigg\{ \tau > n,\ \ \sum_{k=1}^{J_0^n-1}
\textrm{length}({\cal Z}_k^n)\ge \sum_{k=1}^{J_0^n-1} V_k-c n
\bigg\}, $$
since on $E_n$, by Step $(3,n)$ simulator number $0$ would have
used all the bits $\epsilon_j(m)$ for $1\le j\le J_0^n-1$, $1\le
m\le V_j$, which were left unused by the simulators numbered
$1,\ldots,J_0^n-1$, and there cannot be more than $cn$ such.

\begin{sloppy}
Next, observe that on the event $\{\tau>n\}$, because
$(J_0^n,M_0^n)$ $=\textrm{NEXT}^n(J_0^0,M_0^0)$ (the $n$-th iteration of
$\textrm{NEXT}$), we have the equality
$$ J_0^n = \inf\bigg\{ j\ge 0: \sum_{i=0}^j V_i > n \bigg\} =:
\theta_n $$ \end{sloppy}
(we denote the quantity on the right-hand side by $\theta_n$). So we
have shown
$$ E_n \subseteq \bigg\{ \tau > n,\ \ \sum_{k=1}^{\theta_n-1}
\textrm{length}({\cal Z}_k^n)\ge \sum_{k=1}^{\theta_n-1} V_k-c n
\bigg\}. $$
The value of $\theta_n$ is, with probability exponentially close to
$1$, close to $({\bf E}(V_1))^{-1}n$. More precisely, for any $\delta>0$
\begin{eqnarray*}
&\bigg\{ \theta_n > (({\bf E}(V_1))^{-1}+\delta)n \bigg\}
\subseteq \bigg\{ \sum_{0\le i\le (({\bf E}(V_1))^{-1}+\delta)n} V_i \le
n \bigg\}& \\ 
&=
\bigg\{ \sum_{0\le i\le (({\bf E}(V_1))^{-1}+\delta)n} (V_i-{\bf E}(V_1))
\le -{\bf E}(V_1)\cdot \delta n \bigg\}&
\end{eqnarray*}
has probability which decays exponentially by Lemma \ref{associated}
and Lemma \ref{expotails}(iii), and similarly, the probability of
$$ \bigg\{ \theta_n < (({\bf E}(V_1))^{-1}-\delta)n \bigg\}
\subseteq \bigcup_{0\le j< (({\bf E}(V_1))^{-1}-\delta)n}
 \bigg\{ \sum_{0=1}^j V_i > n \bigg\} $$
decays exponentially. So, summarizing the latest developments, we can
now write
\begin{eqnarray*}
E_n &\subseteq& \bigg\{ |\theta_n-({\bf E}(V_1))^{-1}n| > \delta n
\bigg\} \\ & & \cup
\bigcup_{|j-({\bf E}(V_1))^{-1}n|\le \delta n}
\bigg\{ \sum_{k=1}^j \textrm{length}({\cal Z}_k^n) \ge
\sum_{k=1}^j V_k - cn \bigg\}.
\end{eqnarray*}
The first event has probability which decays exponentially by the remarks
above. The second event is a union over linearly many events, each of whose
probability we can only increase by replacing $\textrm{length}({\cal
   Z}_k^n)$ by $T_k'$, an independent copy of $T_{\lambda_1}$,
using Lemma \ref{induction}(iv). So it is enough to prove that
$$
\max_{|j-({\bf E}(V_1))^{-1}n|\le \delta n} {\bf P}\bigg( \sum_{k=1}^j
T_k' \ge \sum_{k=1}^j V_k - cn \bigg) $$
decays exponentially in $n$. For
this, invoking Lemma \ref{tvalue}, choose $c$ and $\delta$ sufficiently
small (the following choices will work: $\delta=({\bf E}(V_1))^{-1}/2$,
$c=({\bf E}(V_1)-{\bf E}(T_{\lambda_1}))/(3(1+{\bf E}(V_1)))$ so that the
inclusion
\begin{multline*}
\bigg\{ \sum_{k=1}^j T_k' \ge
\sum_{k=1}^j V_k - cn \bigg\} \ \subseteq\ 
\bigg\{ \sum_{k=1}^j (T_k'-{\bf E}(T_{\lambda_1})) \ge cj \bigg\} \\ \cup
\bigg\{ \sum_{k=1}^j (V_k-{\bf E}(V_1)) \le -cj \bigg\},
\qquad |j-({\bf E}(V_1))^{-1}n| \le \delta n,
\end{multline*}
will hold, and use Lemma \ref{expotails}(iii).
\end{pf}

\begin{lemma} \label{expo}
$\varphi$ is a finitary homomorphism from $B(p)$ to $B(q)$ with
exponential tails.
\end{lemma}

\begin{pf} From Lemmas \ref{phidist} and \ref{equivariant} it follows
that $\varphi$ is a homomorphism from $B(p)$ to $B(q)$. From the
definition $(\varphi(x))_0 = \beta_{K(0)}(-R_{K(0)}) = \beta_0(-R_0)$
(since $K(0)=0$) we see that
$(\varphi(x))_0$ is determined from the input symbols $(x_i)_{R_0 \le
i\le R_{J_0^\infty+1}+t}$. This proves that $\varphi$ is
finitary, with a coding length $N_\varphi$ that satisfies
\begin{eqnarray*}
N_\varphi &\le& \max(-R_0, R_{J_0^\infty+1}+t) \le
R_{J_0^\infty+2} - R_0 \\ &=&
\sum_{j=1}^{J_0^\infty+1} (R_{j+1}-R_j) + (R_1-R_0).
\end{eqnarray*}
Since by Lemma \ref{wordsindep}, $(R_{j+1}-R_j)_{j>0}$ are i.i.d. with
exponential tails, and $R_1-R_0$ has exponential tails, it follows from
Lemma \ref{expotails}(i),(iv) and Lemma \ref{jkinfty} that
$N_\varphi$ has exponential tails.
\end{pf}

\section{Extension to Markov chains}

We indicate here the changes in the ideas presented above required to prove
Theorem \ref{markov}.  We omit the details of the proofs, which are 
similar to those above. 

\subsection{Coding from a Markov source}

If the Bernoulli source $B(p)$ is replaced by a Markov source, two changes
to the construction are needed. First, the marker pattern ``$211\ldots 1$''
must be replaced with a sequence which is assigned positive probability by
the Markov chain. Here, we must give up the source-universality property,
or make the rather restrictive assumption that all the entries in the
matrix $\alpha$ are positive.

Second, it is necessary to replace the function $F_{n,t}$ that produces
independent unbiased bits from a pattern-free Bernoulli block with a new
function, $F_{n,t}'$ designed to do the same for a pattern-free Markov
block conditioned to begin with a ``1'' (which we assume without loss of
generality to be the last symbol in the left marker) and to end with a
``2'' (the first symbol in the right marker).  Note that the construction
of $F_{n,t}$ used only the symmetries of the distribution
$\tilde{p}_{n,t}$, namely the fact that the space $E_{n,t}$ can be
partitioned into classes of equiprobable elements, since the probability of
an element only depends on the count of the different ${\bf A}$-symbols,
and not on the order of their appearance.

For a Markov source, the probability of an element in $E_{n,t}$ will depend
on the count of adjacent \emph{pairs} of symbols. Thus, there will again be
classes of equiprobable elements. The number of classes, which bounds the
range of the complementary function $G_{n,t}'$ (and hence the amount of
lost entropy -- see the proof of Lemma \ref{tvalue}), is at most
$n^{a^2}$. All the proofs carry through identically to the Bernoulli case.

\subsection{Coding to a Markov chain}

Things get more complicated when the target process is Markov. Here, it is
not enough for each block to independently generate the ${\bf B}$-symbols
required to fill its spaces, since one must make sure that there are the
correct transition probabilities on the boundaries between blocks. We
propose the following solution to this problem.

First, we need a process version of Theorem \ref{simu1}. That is, given a
sequence of random variables $Y_1, Y_2, \ldots$, not necessarily
independent, taking values in a finite alphabet ${\bf B}$, one may
construct a sequence of simulations $(T_k, S_k)_{k\in\nat}$, such that the
stopping times $T_k$ are increasing, and for each $k$, $(T_k, S_k)$ is a
simulation of the distribution of $(Y_1, Y_2, \ldots, Y_k)$ from
independent unbiased bits. Each simulation is efficient, in the sense of
Theorem \ref{simu1}(ii). The proof is obvious, by successively refining the
partition $0=Q_0<Q_1<\ldots<Q_b=1$ of $(0,1)$ used in the proof of
Theorem \ref{simu1}.

Since the Markov matrix $\beta$ is irreducible and aperiodic, there exists
an $m\ge 1$ for which all the entries of $\beta^m$ are positive. For each
block $k\in\z$, the $k$-th simulator will start by generating, using the
independent unbiased bits $U_k$ at her disposal (and, if necessary, bits
from the blocks to her right), $b$ separate Markov trajectories of length
$m$, $(z_{k,i,j})_{1\le i\le b,\ 1\le j\le m}$, such that for each $1\le
i\le b$, $(z_{k,i,j})_{1\le j\le m}$ is a finite Markov chain that starts
with the symbol $i$ and has transition probabilities given by the matrix
$\beta$. We call these finite chains the {\dof preamble chains}. One of
them will later be chosen to fill the first $m$ places in
the $k$-th block, but at the moment we don't know which.

Next, we want the $k$-th simulator to compute symbols to fill the remainder
of her block. This can be done if in the sequences $(z_{k,i,j})$, {\dof
   coupling} was achieved, i.e., if all the symbols $(z_{k,i,m})_{1\le i\le
   b}$ are identical; since in this case, no matter which of the preamble
chains we later choose, we will need the same conditional distribution to
compute the $(m+1)$-th symbol in the block, then the $(m+2)$-th, and so
on.

Because of our choice of $m$, we know that coupling is achieved with a
positive probability. So a positive proportion of the simulators will be
able to continue and compute symbols to fill their blocks, using the nested
simulations described above. But each simulator that filled her block, also
determined for the block to her right which of the preamble chains to use -
the one that corresponds to the last symbol in the block that was filled.
For each block for which the preamble was determined in this second round
of computations, one may now proceed to simulate the remainder of the
block. This then determines the preamble of more blocks, for which the
block remainder is then computed. By iterating this process, the choice of
preamble chains is propagated to the right until a single output sequence
is determined.

The computation of the preamble chains uses a fixed amount of entropy per
block. Since the blocks may be made arbitrarily large in expectation by
choosing a long enough marker, the loss in entropy can be made negligible.
The computation of the remainder of each block uses nested simulations,
which are efficient. Therefore it can be shown that the total entropy loss
is small, and the simulation will terminate a.s. The output sequence is
clearly a stationary process, and by the construction its transition
probabilities are exactly given by the transition matrix $\beta$.
Therefore, it is the stationary Markov shift ${\cal M}(\beta)$. The
construction is a finitary homomorphism, which can be shown to have
exponential tails using large-deviations estimates similarly to the
Bernoulli case.

\section*{Open Problems}

\begin{mylist}
\item
Do there exist source-universal finitary {\em isomorphisms}?  More
specifically, if the Bernoulli sources $B(p),B(p'),B(q)$ all have
equal entropy, does there exist a finitary map (as defined in the
introduction) which is
simultaneously a finitary isomorphism from $B(p)$ to $B(q)$, and from
$B(p')$ to $B(q)$?

{\bf Remark.} If the function is not required to be a finitary map,
the answer to the above question is positive, for a trivial reason.
The function can simply use the
law of large numbers to discern whether the input
is in the almost-sure set
of $B(p)$ or of $B(p')$, and apply one of two
Keane-Smorodinsky \cite{ks:iso}
finitary isomorphisms accordingly.

\item
Construct finitary isomorphisms between general Bernoulli sources
with explicit bounds on the tails. 
(See \cite{harveyperes},\cite{kalikowweiss},\cite{meshalkin}
for such constructions in specific cases).
\end{mylist}

\bigskip

\bigskip \noindent
Nate Harvey \\
Department of Mathematics \\
UC Berkeley \\
CA 94720, USA

\bigskip \noindent
Alexander E. Holroyd \\
Department of Mathematics \\
University of British Columbia \\
Vancouver, BC V6T 1Z2, Canada \\
\texttt{holroyd@math.ubc.ca}

\bigskip \noindent
Yuval Peres \\
Departments of Statistics and Mathematics \\
UC Berkeley \\
CA 94720, USA \\
\texttt{peres@stat.berkeley.edu}

\bigskip \noindent
Dan Romik \\
The Mathematical Sciences Research Institute \\
17 Gauss Way \\
Berkeley, CA 94720-5070 \\
USA \\
\texttt{dromik@msri.org}

\end{document}